\documentclass[proceedings]{stacs}
\stacsheading{2010}{501-512}{Nancy, France}
\firstpageno{501}

\usepackage{amssymb}

\usepackage{xspace}

\usepackage[all,dvips]{xy}

\usepackage{pst-plot, pstricks} 

\newcommand*{\N}[0]{\ensuremath{\mathbb{N}}}

\newcommand*{\domain}[0]{\ensuremath{\text{dom}}}

\newcommand*{\Push}[1]{\ensuremath{\mathrm{push}_{#1}}}
\newcommand*{\Pop}[1]{\ensuremath{{\mathrm{pop}_{#1}}}}
\newcommand*{\Clone}[1]{\ensuremath{{\mathrm{clone}_{#1}}}}
\newcommand*{\Collapse}[0]{\ensuremath{\mathrm{collapse}}}
\newcommand*{\TOP}[1]{\ensuremath{\mathrm{top}_{#1}}}
\newcommand*{\Sym}[0]{\ensuremath{\mathrm{Sym}}}
\newcommand*{\Lvl}[0]{\ensuremath{\mathrm{CLvl}}}
\newcommand*{\Lnk}[0]{\ensuremath{\mathrm{CLnk}}}
\newcommand*{\Op}[0]{\mathrm{OP}}
\newcommand*{\op}[0]{\mathrm{op}}

\newcommand*{\Stacks}[0]{\mathrm{Stck}}

\newcommand*{\Graph}[0]{\mathrm{CPG}}

\newcommand*{\length}[0]{\mathrm{ln}}

\newcommand*{\Loops}[0]{\mathrm{Loops}}

\newcommand*{\Decode}[0]{\mathrm{Dec}}
\newcommand*{\Encode}[0]{\mathrm{Enc}}

\newcommand*{\EncTrees}[0]{\ensuremath{\mathbb{T}_{\mathrm{Enc}}}}

\newcommand*{\LeftTree}[1]{\ensuremath{LT({#1})}}

\newcommand*{\LeftStack}[0]{\ensuremath{\mathrm{LStck}}}

\newcommand*{\T}[0]{\ensuremath{\mathbb{T}}}
\newcommand*{\inducedTreeof}[2]{\ensuremath{{#2}_{#1}}}

\newcommand*{\CPS}[0]{\ensuremath{\mathrm{CPS}}\xspace}
\newcommand*{\CPG}[0]{\ensuremath{\mathrm{CPG}}\xspace}

\newcommand*{\MSO}{\ensuremath{\mathrm{MSO}}\xspace}
\newcommand*{\FO}{\ensuremath{\mathrm{FO}}\xspace}

\newcommand*{\Reach}{\ensuremath{\mathrm{Reach}}}

\newcommand*{\coloneqq}[0]{\ensuremath{\mathrel{\mathop:}=}}

\newcommand*{\trans}[1]{\ensuremath{\mathrel{{\vdash^{#1}}}}}

\newcommand*{\Milestones}[0]{\ensuremath{\mathrm{MS}}}

\newcommand*{\treeLR}[3]{\ensuremath{#1({#2},{#3})}}
\newcommand*{\treeL}[2]{\ensuremath{#1({#2},\emptyset)}}
\newcommand*{\treeR}[2]{\ensuremath{#1(\emptyset,{#2})}}

\begin{document}

\title[Tree-Automaticity of $2$-CPG]{Collapsible Pushdown Graphs of Level
  $\mathbf{2}$ are Tree-Automatic}

\author[lab1]{A. Kartzow}{Alexander Kartzow}
\address[lab1]{TU Darmstadt, Fachbereich Mathematik,
Schlossgartenstr. 7, 64289 
  Darmstadt, Germany}



\begin{abstract}
  \noindent  We show that graphs generated by collapsible pushdown systems of
  level $2$ are tree-automatic. Even when we allow
  $\varepsilon$-contractions and add a reachability
  predicate (with regular constraints) for pairs of configurations,
  the structures remain 
  tree-automatic. Hence, their \FO theories are decidable, even when
  expanded by a reachability predicate.
  As a corollary, we obtain the tree-automaticity of the second level
  of the Caucal-hierarchy. 
\end{abstract}

\maketitle
\section{Introduction}

Higher-order pushdown systems were first introduced by Maslov
\cite{Maslov74,Maslov76} as accepting devices for word languages. 
Later, Knapik et al. \cite{KNU02} studied them as generators for
trees. They obtained an equi-expressivity result for higher-order
pushdown systems and for higher-order recursion schemes that satisfy
the constraint of \emph{safety}, which is a rather unnatural syntactic
condition. Recently, 
Hague et al. \cite{Hague2008} introduced collapsible pushdown systems
as extensions of higher-order pushdown systems and proved that these
have exactly the same power as higher-order recursion schemes as
methods for generating trees. 

Both -- higher-order and collapsible pushdown systems -- also form
interesting devices for generating graphs. 
Carayol and W\"ohrle \cite{cawo03} showed that the graphs generated by
higher-order pushdown systems\footnote{The graph
  generated by a higher-order pushdown system is the
  $\varepsilon$-closure
  of its reachable configurations.}
of level $l$  coincide with the graphs in 
the $l$-th level of the Caucal-hierarchy, a class of graphs
introduced by Caucal \cite{Caucal02}. Every level of this hierarchy is
obtained from the preceding level by applying graph unfoldings
and \MSO interpretations. Both operations preserve the
decidability of the \MSO theory whence the Caucal-hierarchy forms a
rather large class of graphs with decidable \MSO theories. 
If we use collapsible pushdown systems as generators for graphs we
obtain a different situation. Hague et al. showed that even the second
level of the hierarchy contains a graph with undecidable
\MSO theory. But they showed the decidability of the
modal $\mu$-calculus theories of
all graphs in the hierarchy. This turns graphs generated by
collapsible pushdown systems into an interesting class from a model
theoretic point of view. There are few natural classes that share
these properties. In fact, the author only knows one further example,
viz.\ nested pushdown trees. Alur et al.\cite{Alur06languagesof}
introduced
these graphs for $\mu$-calculus model checking purposes. 
We proved in \cite{Kartzow09} that nested pushdown trees also have
decidable 
first-order theories. We gave an effective model checking algorithm
using pumping techniques, but we also proved that nested pushdown
trees are tree-automatic structures. 
Tree-automatic structures were introduced by Blumensath
\cite{Blumensath1999}. 
These structures enjoy decidable first-order theories due to the good
closure properties of finite automata on trees. 

In this paper, we are going to extend our previous result to the
second level of the collapsible pushdown hierarchy. All graphs of the
second level are tree-automatic. This subsumes our previous result as
nested pushdown trees are first-order interpretable in collapsible
pushdown graphs of level two. 
Furthermore, we show that collapsible pushdown graphs of level $2$
are still tree-automatic when
expanded by a reachability predicate, i.e., by the binary relation
which contains all pairs of configurations such that there is a path
from the first to the second configuration. Thus, first-order logic
extended by reachability
predicates is decidable on level $2$ collapsible pushdown graphs. 

In the next section, we introduce the necessary notions concerning
tree-automaticity and in Section \ref{STACS:SecCPG} we define
collapsible pushdown graphs. We explain the translation of
configurations into trees in Section
\ref{STACS:SecEncoding}. 
Section \ref{STACS:SecCertificates} 
is a sketch of the proof that this translation yields tree-automatic
representations of collapsible pushdown graphs, even when enriched with
certain regular reachability predicates.
The last section contains some concluding remarks about questions
arising from our result.

\section{Preliminaries}
\noindent
We write \MSO for monadic second order logic and \FO for first-order logic.
\noindent
For words $w_1,w_2\in \Sigma^*$, we write $w_1 \sqcap w_2$ for the
greatest common prefix of $w_1$ and $w_2$.
\noindent
A \emph{$\Sigma$-labelled tree} is a function $T:D\rightarrow \Sigma$
for a finite $D\subseteq\{0,1\}^*$ which is closed under prefixes.

\noindent
For $d\in D$ we denote by
$\inducedTreeof{d}{T}$ the \emph{subtree rooted at $d$}.

Sometimes it is useful to define trees inductively by describing their
left and right subtrees. For this purpose we fix the following notation.
Let $\hat T_0$ and $\hat T_1$ be $\Sigma$-labelled trees and
$\sigma\in\Sigma$. Then we write 
$T\coloneqq \treeLR{\sigma}{\hat T_0}{ \hat T_1}$
for the $\Sigma$-labelled tree $T$ with the following three
properties
\begin{align*}
  &1.\ T(\varepsilon) = \sigma,& &2.\ \inducedTreeof{0}{T} = \hat
  T_0 \text{, and }& &3.\ \inducedTreeof{1}{T} = \hat T_1\enspace.&
\end{align*}

In the rest of this section, we briefly present the notion of a
tree-automatic structure as introduced by Blumensath \cite{Blumensath1999}.

The \emph{convolution} of two
$\Sigma$-labelled trees $T$ and $T'$ is given by a function
\begin{align*}
  T\otimes T' : \domain(T)\cup\domain(T') \rightarrow
  (\Sigma \cup \{\Box\} )^2
\end{align*}
where $\Box$ is a new symbol for padding and
\begin{align*}
  (T\otimes T')(d) \coloneqq
  \begin{cases}
    (T(d),T'(d)) & \text{ if } d\in \domain(T)\cap
    \domain(T') \\ 
    (T(d), \Box) & \text{ if }d\in \domain(T)\setminus
    \domain(T') \\ 
    (\Box, T'(d)) & \text{ if }d\in \domain(T') \setminus
    \domain(T)
  \end{cases}
\end{align*}
By ``tree-automata'' we mean a nondeterministic finite automaton that
labels a finite tree top-down.

\begin{definition}
  A structure $\mathfrak{B} = (B,E_1, E_2, \ldots, E_n)$ with domain
  $B$ and binary relations
  $E_i$ is \emph{tree-automatic} if there are tree-automata $A_B,
  A_{E_1}, A_{E_2}, 
  \ldots, A_{E_n}$ and a bijection \mbox{$f: L\rightarrow B$} 
  for $L$ the language accepted by
  $A_B$ such that the
  following hold.
  For $T, T'\in L$, the automaton $A_{E_i}$ accepts
    $T\otimes T'$ if and only if \mbox{$\big(f(T), f(T')\big)\in E_i$.} 
\end{definition}
Tree-automatic structures form a nice class because automata theoretic
techniques may be used to decide first-order formulas on these
structures:

\begin{lemma}[\cite{Blumensath1999}]
  If $B$ is tree-automatic, then its first-order theory
  is decidable.
\end{lemma}

We will use the classical result that regular sets of trees are
$\MSO$ definable.

\begin{theorem}[\cite{ThatcherW68}, \cite{Doner70}]
  For a set $\T$ of finite $\Sigma$-labelled trees,
  there is a tree 
  automaton recognising $\T$ if
  and only if $\T$ is \MSO definable.
\end{theorem}

\section{Definition of Collapsible Pushdown Graphs (CPG)}
\label{STACS:SecCPG}
In this section we define our notation of collapsible pushdown systems. For
a more comprehensive introduction, we refer the 
reader to \cite{Hague2008}.

\subsection{Collapsible Pushdown Stacks}

First, we provide some terminology concerning stacks of (collapsible)
higher-order pushdown systems.
We write $\Sigma^{*2}$ for $(\Sigma^*)^*$ and $\Sigma^{+2}$ for
$(\Sigma^+)^+$. We call an  
$s\in\Sigma^{*2}$ a $2$-word. 

Let us fix a $2$-word $s\in \Sigma^{*2}$ which
consists of an ordered list 
\mbox{$w_1, w_2, \ldots, w_m\in\Sigma^*$}. 
We
separate the words of this list by colons writing 
$s=w_1: w_2 : \ldots: w_m$.
By $\lvert  s \rvert$ we denote the number of words $s$
consists of, i.e., $\lvert  s\rvert = m$.

For another word $s' = w_1':w_2':\ldots :w'_n \in \Sigma^{*2}$, we write 
$s:s'$ for the concatenation  $w_1: w_2 : \ldots :w_m:w_1':w_2':\ldots:w_n'$.

If $w\in\Sigma^{*}$, we write $[w]$ for the $2$-word
that consists of a list of one word which is 
$w$.

A level $2$ collapsible pushdown stack is a special element of
\mbox{$(\Sigma\times\{1, 2\}\times\N)^{+2}$} that is generated by
certain stack operations 
from an initial stack which we introduce in the following definitions.
The natural numbers following the stack symbol represent the so-called
\emph{collapse pointer}: every element in a collapsible pushdown stack has a
pointer to some substack and applying the collapse operation
returns the substack to which the topmost symbol of the stack
points. Here, the first number denotes the \emph{collapse level}.  
If it is $1$ the collapse pointer always points to the symbol below
the topmost symbol and the collapse operations just removes the
topmost symbol. 
The more interesting case is when the collapse level of the topmost
symbol of the stack $s$ is $2$. Then the
stack obtained by the collapse contains the first $n$ words of $s$
where $n$ is the second number in the topmost element of $s$. 

The initial level $1$ stack is $\bot_1\coloneqq(\bot,1,0)$ and the
initial level $2$ stack is  $\bot_2\coloneqq[\bot_1]$. 

For $k\in\{1,2\}$ and for a
$2$-word $s=w_1:w_2:\ldots: w_n\in(\Sigma\times\{1, 2\} 
\times\N)^{+2}$ such that $w_n=a_1 a_2 \ldots a_m$ with
$a_i\in\Sigma\times\{1,2\}\times\N$ for all $1\leq i\leq m$:
\begin{itemize}
\item we define the \emph{topmost $(k-1)$-word of $s$} as $\TOP{k}(s)\coloneqq
  \begin{cases}
    w_n &\text{if } k=2\\
    a_m &\text{if } k=1
  \end{cases}$
\item for $\TOP{1}(s)=(\sigma,i,j)\in \Sigma\times\{1, 2\}\times\N$,  
  we define the \emph{topmost symbol}
  \mbox{$\Sym(s)\coloneqq\sigma$}, the  \emph{collapse-level of the
    topmost element}
  $\Lvl(s)\coloneqq i$, and the \emph{collapse-link of the topmost element}
  $\Lnk(s)\coloneqq j$. 
\end{itemize}
For $s$, $w_n$ and $k$ as before,
$\sigma\in\Sigma\setminus\{\bot\}$, and $w_n':=a_1 \ldots a_{m-1}$, we
define the stack operations 
\begin{align*}
  \Pop{k}(s)\coloneqq
  &\begin{cases}
    w_1: w_2 : \ldots : w_{n-1} & \text{if } k=2, n\geq2\\
    w_1: w_2 : \ldots : w_{n-1}: w_n' & \text{if } k=1, m\geq 2\\
    \text{undefined} & \text{otherwise} 
  \end{cases}  \\
  \Clone{2}(s)\coloneqq&\ w_1: w_2 : \ldots : w_{n-1}: w_n : w_n \\
  \Push{\sigma,k}(s)\coloneqq
  &\begin{cases}
    w_1 : w_2 : \ldots :w_n(\sigma,2,n-1) &\text{ if k=2}\\
    w_1 : w_2 : \ldots :w_n(\sigma,1,m) &\text{ if k=1}
  \end{cases} \\
  \Collapse{}(s)\coloneqq 
  &\begin{cases}
    w_1: w_2 : \ldots : w_r & \text{if } \Lvl(s)=2, \Lnk(s) = r > 0\\
    \Pop{1}(s) & \text{if } \Lvl(s)=1\\
    \text{undefined} & \text{otherwise} 
  \end{cases}    
\end{align*}
The \emph{set of level $2$-operations} is
$\Op\coloneqq\left\{\Push{\sigma,1},\Push{\sigma,2}, 
\Clone{2}, \Pop{1}, \Pop{2},
\Collapse{}\right\}$.   
The \emph{set of level $2$ stacks}, $\Stacks(\Sigma)$, is the smallest set
that contains $\bot_2$ and is closed under all operations from $\Op$.

Note that $\Collapse$- and $\Pop{k}$-operations are only allowed if
the resulting stack is in $(\Sigma^+)^+$. This avoids the special
treatment of empty words or stacks. 
Furthermore, a $\Collapse$ on level $2$ summarises a non-empty
sequence of $\Pop{2}$-operations. For example, starting from $\bot_2$,
we can apply a $\Clone{2}$, a $\Push{\sigma,2}$, a
$\Clone{2}$, and finally a $\Collapse$. This sequence first creates a level $2$
stack that contains $3$ words and then performs the collapse and ends in the
initial stack again. This example shows that
$\Clone{2}$-operations are responsible for the fact that 
collapse-operations on level $2$ may remove more than one word from the
stack.   
 
For $s,s'\in\Stacks(\Sigma)$, we call $s'$ a substack of $s$ if
there are  $n_1, n_2\in\N$ such that
\mbox{$s' = \Pop{1}^{n_1}( \Pop{2}^{n_2}(s))$.}
We write $s' \leq s$ if $s'$ is a substack of $s$.

\subsection{Collapsible Pushdown Systems and Collapsible Pushdown Graphs}

Now we introduce
collapsible pushdown systems and graphs (of level $2$) which are
analogues of pushdown systems and pushdown graphs
using collapsible pushdown stacks instead of ordinary stacks. 

\begin{definition}
  A \emph{collapsible pushdown system} of level $2$ (\CPS) is
  a tuple 
  \mbox{$S = (\Sigma, Q, \Delta, q_0)$} where $\Sigma$ is a finite stack
  alphabet with $\bot\in\Sigma$, $Q$ a
  finite set of states, $q_0\in Q$ the initial state, and $\Delta\subseteq
  Q\times \Sigma \times Q \times \Op$ the transition relation.

  For $q\in Q$ and $s\in\Stacks(\Sigma)$  the  pair $(q,s)$ is called 
  a \emph{configuration}. 
  We 
  define labelled
  transitions on pairs of configurations by setting $(q_1,s)
  \trans{(q_2,op)} (q_2, t)$ if there is a
  $(q_1, \sigma, q_2, op)\in\Delta$ such that $\Sym(s) = \sigma$ and
  $op(s)=t$. The union of the labelled transition
  relations is denoted as 
  \mbox{$\trans{} \coloneqq\bigcup_{l\in Q\times \Op}
  \trans{l}$.}
  We set $C(S)$ to be the set of all
  configurations that are reachable from $(q_0,\bot_2)$ via
  $\trans{}$-paths. We call $C(S)$ the set of \emph{reachable}
  or \emph{valid} configurations. The 
  \emph{collapsible pushdown graph (\CPG) generated by $S$} is 
  \begin{align*}
    \Graph(S)\coloneqq\left(C(S), (C(S)^2\cap
    \trans{\ell})_{\ell\in Q\times\Op}\right)
  \end{align*}
\end{definition}

\begin{example}
  The following example of a collapsible pushdown graph of level $2$
  is taken from 
  \cite{Hague2008}. Let $Q\coloneqq\{0,1,2\}, \Sigma\coloneqq \{\bot,a\}$, and
  $\Delta$ given by $(0,*,1,\Clone{2})$,
  $(1,*,0,\Push{a,2})$, $(1,*,2,\Push{a,2})$,
  $(2,a,2,\Pop{1})$, and $(2,a,0,\Collapse)$, where $*$ denotes any
  letter in $\Sigma$. 
  In our picture (see Figure \ref{STACSfig:CPGExample}), the labels
  are abbreviated as follows:
  $\mathrm{cl}\coloneqq (1,\Clone{2})$, 
  \mbox{$a\coloneqq (0,\Push{a,2})$,} 
  $a'\coloneqq (2,\Push{a,2})$, 
  $p\coloneqq (2,\Pop{1})$, and
  $\mathrm{co}\coloneqq (0,\Collapse)$. \\
  \begin{figure}[h]
    \centering
				\resizebox{\textwidth}{!}{
    $
    \begin{xy}
      \xymatrix@R=12pt@C=10pt{
        0, \bot \ar[r]^-{\mathrm{cl}}& 
          1, \bot:\bot \ar[r]^-{ a} \ar[d]^{ a'}& 
          0, \bot:\bot a \ar[r]^-{\mathrm{cl}} &
          1, \bot:\bot a:\bot a  \ar[r]^-{a} \ar[d]^-{a'}& 
          0,\bot:\bot a: \bot aa \ar[r]^-{\mathrm{cl}}  
          & 1, \bot:\bot a : \bot aa:\bot aa \ar[d]^{a'} \ar[r]^-{ a}& \ldots \\
        & 
          2, \bot:\bot a  \ar[d]^p \ar[ul]^{\mathrm{co}}&  
          & 
          2, \bot:\bot a:\bot aa \ar[d]^p \ar[ul]^{\mathrm{co}}&  
          & 
          2, \bot:\bot a : \bot aa:\bot aaa\ar[d]^p
            \ar[ul]^{\mathrm{co}}
          & \ldots \\
        & 
          2, \bot:\bot  &  
          &
          2, \bot:\bot a:\bot a  \ar[d]^p \ar[uulll]^{\mathrm{co}}&  
          & 
          2, \bot:\bot a : \bot aa:\bot aa
          \ar[d]^p\ar[uulll]^(.415){\mathrm{co}}
          & \ldots \\
        &                &  
          & 
          2, \bot:\bot a:\bot  &  
          & 
          2, \bot:\bot a : \bot aa:\bot a \ar[d]^p\ar[uuulllll]^(.3){\mathrm{co}}
          & \ldots \\
        &              
          &  
          &
          &
          &
          2, \bot:\bot a : \bot aa:\bot  & \ldots 
      }
    \end{xy}
    $
				}
    \caption{Example of a collapsible pushdown graph}
    \label{STACSfig:CPGExample}
  \end{figure}
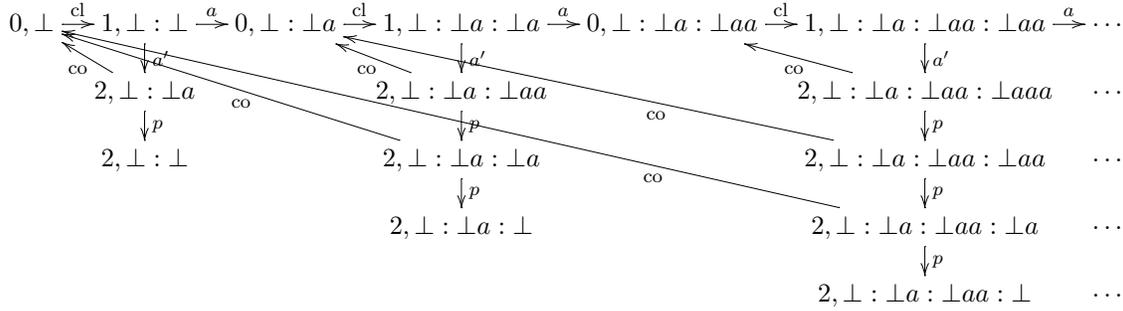
\end{example}

\begin{remark} \label{STACS:LmCPSDecidable}
  Hague et al. \cite{Hague2008} showed that
  modal $\mu$-calculus model checking on level $n$ \CPG is $n$-EXPTIME
  complete. Note that there is
  an $\MSO$ interpretation which turns the graph of the previous
  example into a grid-like
  structure. Hence its \MSO theory is undecidable.
\end{remark}

The next definition introduces runs of collapsible pushdown systems.

\begin{definition}
  Let $S$ be a \CPS.
  A run $r$ of $S$ of length $n$ is a 
  function 
  \begin{displaymath}
       r:\{0, 1, 2, \ldots ,n\} \rightarrow
    Q\times(\Sigma\times\{1, 2\}\times\N)^{*2}
    \text{ such that }
    r(0) \trans{} r(1) \trans{} \cdots \trans{} r(n). 
  \end{displaymath}
  We write $\length(r)\coloneqq n$ and
  call $r$ a run from $r(0)$ to $r(n)$. We say $r$ visits
  a stack $s$ at $i$ if $r(i)=(q,s)$.
  
  For runs $r, r'$ of length $n$ and $m$, respectively, with
  $r(n)=r'(0)$, we define the composition $r\circ r'$ of $r$ and $r'$
  in the obvious manner.
\end{definition}

\begin{remark}
  Note that we do not require runs to start in the initial configuration.
\end{remark}

\section{Encoding of Collapsible Pushdown Graphs in Trees}
\label{STACS:SecEncoding}
\begin{figure}
  \centering
  $
  \begin{xy}
    \xymatrix@=0.2mm{
      & f \\
      & e & g & & i\\
      b& d & d & d & h & & j & l\\
      a & c & c & c & c & c & c & k\\
      \bot & \bot & \bot & \bot & \bot & \bot & \bot & \bot
      \save "1,2"."4,4"*[F.]\frm{}
      \save "2,5"."4,5"*[F.]\frm{}
      \save "4,6"."4,6"*[F.]\frm{}
      \save "3,7"."4,7"*[F.]\frm{}
    }
  \end{xy}
  $
  \caption{Example of blocks in a stack. These form a $c$-blockline.}
  \label{STACSfig:Blocks}
\end{figure}
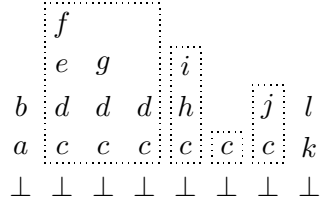

In this section we
prove that \CPG are tree-automatic. 
For this purpose we have to encode stacks in trees. 
The idea is to divide a stack into \emph{blocks} and to
encode different blocks in different subtrees.
The crucial observation is that every stack is a list of words that
share the same first letter. A block is a maximal list of words
in the stack that share the same two first letters\footnote{see Figure
\ref{STACSfig:Blocks} for an example of blocks and Definition
\ref{STACS:DefBlock} for their formal definition}. If we remove the first
letter of every word  of such a block, the resulting $2$-word
decomposes again as a list of blocks. Thus, we can inductively
carry on to decompose parts of a stack into blocks and code every
block in a different subtree. The roots of these subtrees are
labelled with the first letter of the corresponding block. 
This results in a tree in which every initial left-closed
path represents one word of the stack. By
left-closed, we mean that the last element of the path has no left
successor.

It turns out that -- via this encoding -- each stack operation 
corresponds to a simple $\MSO$-definable tree-operation. 
The main difficulty is to provide a
tree-automaton that checks whether there is a run to the configuration
represented by some tree. This problem is addressed in Section
\ref{STACS:SecCertificates}.

As already mentioned, the encoding works by dividing stacks
into blocks.  
The following definition makes our notion of blocks precise. 
For $w\in\Sigma^*$ and $s=w_1:w_2:\ldots:w_n \in\Sigma^{*2}$, we
write $s'\coloneqq w\mathrel\backslash s$ 
for $s'=[ww_1] : [ww_2] : \ldots: [ww_n]$. 

\begin{definition}
  [$\sigma$-block(line)]\label{STACS:DefBlock} For $\sigma\in\Sigma$,
  we call
  $b\in\Sigma^{*2}$ a \emph{$\sigma$-block} if
  $b=[\sigma]$ or $b=\sigma\tau\mathrel{\backslash} s'$ for some
  $\tau\in\Sigma$ and $s'\in\Sigma^{*2}$.
  See Figure \ref{STACSfig:Blocks} for examples of blocks. 
  If $b_1, b_2,\ldots, b_n$ are $\sigma$-blocks, then we call   $b_1:b_2:\ldots :b_n$ a \emph{$\sigma$-blockline}.
\end{definition}

Note that every 
stack in $\Stacks(\Sigma)$ forms a $(\bot,1,0)$-blockline. 
Furthermore, every blockline $l$ decomposes uniquely as
\mbox{$l=b_1: b_2: \ldots: b_n$} of maximal blocks $b_i$ in $l$. 
Another crucial observation is that a $\sigma$-block
$b\in\Sigma^{*2}\setminus \Sigma$
decomposes as $b=\sigma\mathrel{\backslash}l$ for some blockline $l$
and we say $l$ is the induced blockline of $b$. For $b\in\Sigma$ the
induced blockline of $[b]$ is just the empty $2$-word.

Now we encode a $(\sigma,n,m)$-blockline $l$ in a tree by labelling the root
with $(\sigma,n)$, by encoding the blockline induced by the first block of
$l$ in the left subtree, and by encoding the
rest of the 
blockline in the right
subtree. In order to avoid repetitions,
we do not repeat the
symbol $(\sigma,n)$ in the right subtree, but 
replace it by the default letter $\varepsilon$. 

\begin{figure}[t]
  \centering
  $
  \begin{xy}
    \xymatrix@R=0pt@C=0pt{
      &       & (c,2,1) &         & (e,1,3) & \\
      &(b,2,0)& (b,2,0) & (c,1,2) & (d,2,3) & \\
      &(a,2,0)& (a,2,0) & (a,2,2) & (a,2,2) & (a,2,2) \\
      & (\bot,1,0) & (\bot,1,0) &(\bot,1,0) & (\bot,1,0) & (\bot,1,0)
      }
  \end{xy}$
  \hskip 1cm
  $\begin{xy}
    \xymatrix@R=9pt@C=3pt{
           & c,2 &       & e,1      &  \\
      b,2 \ar[r]& \varepsilon\ar[u] & c,1 & d,2\ar[u] & \\
      a,2 \ar[u] &  & a,2\ar[r]\ar[u] & \varepsilon\ar[r]\ar[u]
      &\varepsilon \\
       \bot,1 \ar[rr]\ar[u] &  &\varepsilon\ar[u] &  & &      
      }
  \end{xy}
  $
  \caption{A stack $s$ and its encoding $\Encode(s)$: right
    arrows lead to $1$-successors (right successors), upward arrows
    lead to $0$-successors (left successors).}
  \label{STACSfig:Encoding}
\end{figure}
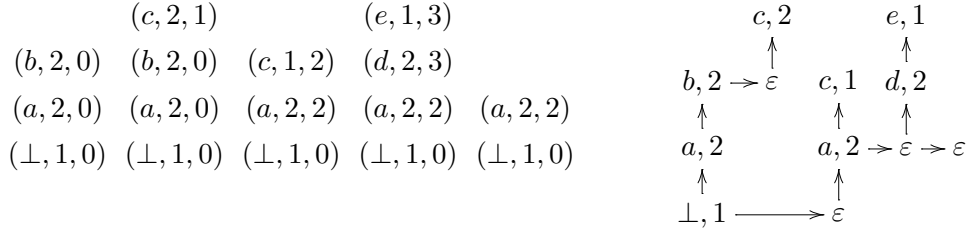

\begin{definition} 
  Let $s = w_1:w_2:\ldots:w_n\in(\Sigma\times\{1,2\}\times\N)^{+2}$ be a
  $(\sigma,l,k)$-blockline. Let $w_i'$ be words
  such that 
  \mbox{$s= (\sigma,l,k) \mathrel\backslash [w_1': w_2' :\ldots
    :w_n']$} and set $s'\coloneqq w_1': w_2' : \ldots : w_n'$. As an
  abbreviation we write  
  $_hs_i\coloneqq w_h:w_{h+1}:\ldots:w_i$.
  Furthermore, let  
  $w_1:w_2:\ldots: w_j$ be a maximal block of $s$. Note that $j>1$ implies 
  $w_{j'}=(\sigma,l,k)(\sigma',l',k') w_{j'}''$ for all $j'\leq j$,
  some fixed
  $(\sigma',l',k')\in \Sigma\times\{1,2\}\times\N$, and appropriate
  $w_{j'}''\in\Sigma^*$. 
  For 
  $\rho\in\big(\Sigma\times\{1,2\}\big)\cup\{\varepsilon\}$, we define
  recursively the
  $\big(\Sigma\times\{1,2\}\big)\cup\{\varepsilon\}$-labelled
  tree $\Encode(s,\rho)$ via
  \begin{align*}
    \Encode(s,\rho)\coloneqq
    \begin{cases}
      \rho & \text{if } \lvert w_1\rvert=1, n=1\\
      \treeR{\rho}{\Encode( _2s_n,\varepsilon)}
      &\text{if } \lvert w_1 \rvert =1, n>1\\
      \treeL{\rho}{\Encode( _1s_n',(\sigma',l'))}
      &\text{if } j=n, \lvert w_1\rvert >1  \\
      \treeLR{\rho}{\Encode( _1s_j',(\sigma',l'))}{\Encode(
        _{j+1}s_n,\varepsilon)}      
      &\text{otherwise.}
    \end{cases}
  \end{align*}
  $\Encode(s)\coloneqq \Encode(s,(\bot,1))$ is called the (tree-)encoding of
  the stack $s\in\Stacks(\Sigma)$.
\end{definition}
Figure \ref{STACSfig:Encoding} shows a configuration and its encoding. 

\begin{remark}\label{STACS:milestonesInEncoding}
  In this encoding, the first block of a $(\sigma,l,k)$-blockline is
  encoded in a subtree
  whose root $d$ is labelled $(\sigma,l)$. 
  We can restore
  $k$ from the position of $d$ in the tree $\Encode(s)$ as follows.
  If $l=1$ then $k=\lvert d \rvert_0$, i.e., the number of occurrences
  of $0$ in $d$. This is due to the fact that
  level $1$ links always point to the preceding letter and that we
  always introduce a left-successor tree in order to encode letters
  that are higher in the stack. 
  
  The case $l=2$ needs some closer inspection.
  Assume that some $d\in T:=\Encode(s)$ is labelled $(\sigma,2)$. Then it
  encodes a letter $(\sigma,2,k)$ and this is not a cloned
  element. Thus, $k$ equals the numbers of words to the left of this letter
  $(\sigma,2,k)$. We claim that
  $k=\left\lvert \big\{e\in T\cap\{0,1\}^*1: e\leq_{lex} d
  \big\}\right\rvert$.
  The existence of a pair $e,e1\in T$ corresponds to the fact
  that there is some blockline consisting of blocks
  $b_1:b_2:\ldots:b_n$ with $n\geq 2$ such that $b_1$ is encoded in
  $\inducedTreeof{e}{T}\setminus\inducedTreeof{e1}{T}$ and
  $b_2:\ldots: b_n$ is encoded in $\inducedTreeof{e1}{T}$. 
  By induction, one easily sees that for each such pair $e,e1\in T$
  all the letters that are in words left of the letter encoded by $e1$
  are encoded in lexicographically smaller elements. 
  Furthermore, the size of $((0^*)1)^*\cap T$ corresponds to the
  number of words in $s$ since the introduction of a $1$-successor
  corresponds to the separation of the first block of some blockline
  from the other blocks. Each of these separation can also be seen as
  the separation of the last word of the first block from the first
  word of the second block of this blockline. Note that we separate
  two words that are next to each other in exactly one blockline. 
  Putting these facts together our claim is proved.
  
  Another view on this correspondence is the 
  bijection
  $f:\{1,2,\ldots,\lvert s \rvert\} \rightarrow R$ where
  \mbox{$R:=((0^*)1)^*\cap
  \domain(T)$} and $i$ is mapped to the $i$-th element of $R$ in
  lexicographic order. $f(i)$ is exactly the position where the
  $(i-1)$-st word is separated from the $i$-th one for all $i\geq 2$. 
  In order to state the properties of $f$, we need some more notation.
  We write $\pi$ for the canonical projection
  $\pi:(\Sigma\times\{1,2\}\times\N)^* \rightarrow
  (\Sigma\times\{1,2\})^*$ and $w_i$ for the $i$-th word of
  $s$. Furthermore, let $w_i'$ be a word such that, $w_i=(w_i\sqcap
  w_{i-1})\circ w_i'$ (here we set $w_0:=\varepsilon$). 
  Then the word along the path\footnote{By the
    word along a path from one node to another we mean the word
    consisting of the non $\varepsilon$-labels along this path.}
  from the root to $f(i)$ is exactly
  $\pi(w_i\sqcap w_{i-1})$ for all $2 \leq i\leq \lvert s \rvert $ and
  the path from $f(j)$ to $f(j)\circ0^m$ for maximal $m\in\N$ is
  $\pi(w_j')$ for all $1\leq j \leq \lvert s \rvert$. 
\end{remark}

In order to encode a configuration $c:=(q,s)$, we
add $q$ as a new root of the tree and attach
the encoding of $s$ as the left subtree, i.e.,
\mbox{$\Encode(c)\coloneqq 
\treeL{q}{\Encode(s)}$}.

The image of this encoding function contains only trees of a very
specific type. We call this class $\EncTrees$. 
In the next definition we state the characterising properties of
$\EncTrees$. This class is  \MSO definable, whence
automata-recognisable. 
\begin{definition} \label{STACS:DefEncodingTrees}
  Let $\EncTrees$ be the class of all trees $T$ that satisfy
  the following conditions.
  \begin{enumerate}
  \item  The root of $T$ is labelled by some element of $Q$
    ($T(\varepsilon)\in Q$).
  \item Every element of the form $\{0,1\}^*0$ is labelled by some
    $(\sigma,l)\in\Sigma\times\{1,2\}$; especially, $T(0)=(\bot,1)$ and
    there are no other occurrences of $(\bot,1)$ or $(\bot,2)$. 
  \item Every element of the form $\{0,1\}^*1$ is labelled by
    $\varepsilon$. 
  \item $1\notin\domain(T)$,   $0\in\domain(T)$.
  \item \label{STACS:fifthofDefEnc} 
    For all $t\in T$, if
    $T(t0)=(\sigma,1)$ then 
    $T(t10)\neq(\sigma,1)$. 
  \end{enumerate}
\end{definition}
\begin{remark}
  Note that (\ref{STACS:fifthofDefEnc}) holds as
  $T(t0)=T(t10)=(\sigma,1)$ would imply that the 
  subtree rooted at $t$ encodes a blockline $l$ such that the first
  block of $l$ induces a $(\sigma,1,n)$-blockline and the second one
  induces a 
  $(\sigma,1,m)$-blockline. But as level $1$ links always point to the
  preceding letter, $n$ and $m$ are equal to the length of the
  prefix of $l$ in the stack plus $1$, i.e., if $T$ encodes a stack
  $s$ then $s=s_1: [w\mathrel{\backslash} l] : s_2$ and $n=m=\lvert w
  \rvert +1$. This would contradict the maximality of the blocks in the
  encoding.
\end{remark}

\begin{remark} \label{STACS:Bijective}
  $\Encode:Q\times \Stacks(\Sigma) \rightarrow \EncTrees$ is a 
  bijection and we denote its inverse by $\Decode$. 
\end{remark}

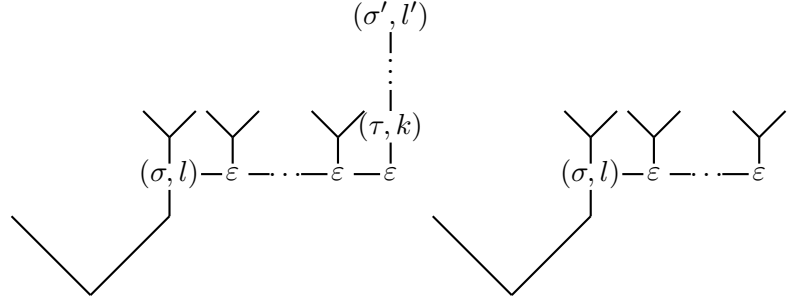
\begin{figure}
  \centering
    \psset{unit=.7mm} 
    \begin{pspicture}(-15,0)(95,50)
      \psline(-15,15)(0,0)
      \psline(15,15)(0,0)
      \psline(15,15)(15,20) \rput(15,23){$(\sigma,l)$}
      \psline(15,25)(15,30)
      \psline(15,30)(10,35)
      \psline(15,30)(20,35)
      \psline(21,23)(25,23) \rput(27,23){$\varepsilon$}
      \psline(27,25)(27,30)
      \psline(27,30)(22,35)
      \psline(27,30)(32,35)
      \psline(30,23)(34,23) \rput(37,23){$\ldots$}
      \psline(40,23)(45,23)\rput(47,23){$\varepsilon$}
      \psline(47,25)(47,30)
      \psline(47,30)(52,35)
      \psline(47,30)(42,35)
      \psline(50,23)(55,23)\rput(57,23){$\varepsilon$}
      \psline(57,25)(57,29)\rput(57,32){$(\tau,k)$}
      \psline(57,35)(57,39)\rput(57,44){$\vdots$}
      \psline(57,46)(57,50)\rput(57,53){$(\sigma',l')$}

      \psline(65,15)(80,0)
      \psline(95,15)(80,0)
      \psline(95,15)(95,20) \rput(95,23){$(\sigma,l)$}
      \psline(95,25)(95,30)
      \psline(95,30)(90,35)
      \psline(95,30)(100,35)
      \psline(101,23)(105,23) \rput(107,23){$\varepsilon$}
      \psline(107,25)(107,30)
      \psline(107,30)(102,35)
      \psline(107,30)(112,35)
      \psline(110,23)(114,23) \rput(117,23){$\ldots$}
      \psline(120,23)(125,23)\rput(127,23){$\varepsilon$}
      \psline(127,25)(127,30)
      \psline(127,30)(132,35)
      \psline(127,30)(122,35)
    \end{pspicture}   
  \caption{$\Pop{2}$-operation}
  \label{STACSfig:Pop2}
\end{figure}

\begin{figure}
  \centering
    \psset{unit=.7mm} 
    \begin{pspicture}(-15,0)(95,40)
      \psline(-15,15)(0,0)
      \psline(10,10)(0,0)
      \psline(10,10)(15,10)
      \psline(15,10)(15,15)
      \rput(15,20){$\vdots$}
      \psline(15,22)(15,27)
      \rput(15,30){$(\sigma,2)$}
      \psline(15,32)(15,37)
      \psline(15,37)(10,42)
      \psline(15,37)(20,42)
      \psline(21,30)(25,30) \rput(27,30){$\varepsilon$}
      \psline(27,32)(27,37)
      \psline(27,37)(22,42)
      \psline(27,37)(32,42)

      \psline(30,30)(34,30) \rput(37,30){$\ldots$}
      \psline(40,30)(45,30)\rput(47,30){$\varepsilon$}
      \psline(47,32)(47,37)
      \psline(47,37)(52,42)
      \psline(47,37)(42,42)
      \psline(50,30)(55,30)\rput(57,30){$\varepsilon$}

      \psline(65,15)(80,0)
      \psline(90,10)(80,0)
    \end{pspicture}   
  \caption{$\Collapse$-operation of level $2$.}
  \label{STACSfig:Col}
\end{figure}
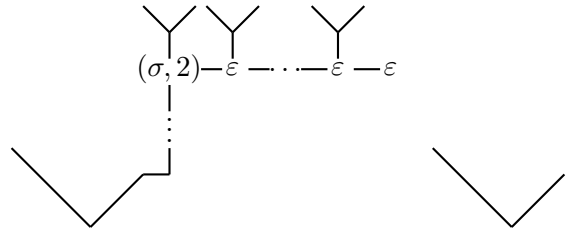

Our encoding turns the transitions of a \CPG into regular
tree-operations. The tree-operations
corresponding to $\Pop{2}$ and $\Collapse$ can be seen in Figures
\ref{STACSfig:Pop2} and  \ref{STACSfig:Col}.
For the $\Pop{2}$, note that if $v_1$ is the $0$-successor of $v_0$
then $v_0$ and $v_1$ encode symbols in the same word
of the encoded stack. As a $\Pop{2}$ removes the rightmost word, we
have to remove all the nodes encoding information about this word. As
the rightmost leaf corresponds to the topmost symbol of the stack, we
have to remove this leaf and all its $0$-ancestors. 

For the $\Collapse$ (on level $2$), we note that each $\varepsilon$
represents a cloned 
element. The $\Collapse$ induced by such an element produces the same
stack as a $\Pop{2}$
of its original version. The original symbol of
the rightmost leaf is its first ancestor not labelled by
$\varepsilon$. 

Note that the operations corresponding to $\Pop{2}$ and
$\Collapse$ are clearly $\MSO$ definable. All other transitions in
$\CPG$ correspond to $\MSO$ definable tree-operations, too. Due to space
restrictions we skip the details.


\begin{lemma}\label{STACS:OperationsRegular}
  Let $C$ be the set of encodings of configurations of a \CPS $S$.
  Then there are automata $A_{(q,\op)}$ for all $q\in Q$
  and all $\op\in\Op$
  such that for all $c_1, c_2\in C$
  \begin{displaymath}
    A_{(q,\op)} \text{ accepts } \Encode(c_1)\otimes \Encode(c_2)
    \text{\quad iff\quad}  c_1 \trans{(q,\op)} c_2\enspace. 
  \end{displaymath}
\end{lemma}

\section{Recognising Reachable Configurations}
\label{STACS:SecCertificates}

We show that $\Encode$ maps the reachable
configurations of a given $\CPS$ to a regular set. 
For this purpose we introduce milestones of a stack
$s$. It turns out that these are exactly those substacks of $s$
that every 
run to $s$ has to visit. Furthermore, the milestones of $s$ are
represented by the nodes of $\Encode(s)$: 
with every $d\in \Encode(s)$, we can associate a subtree of $s$ which
encodes a milestone. Furthermore, the substack relation 
on the milestones corresponds exactly to the lexicographical order
$\leq_{lex}$ of 
the elements of $\Encode(s)$. 
For every \mbox{$d\in \Encode(s)$} we can guess the
state in which the corresponding milestone is visited for the last
time by some run to $s$ and we can check the correctness of this
guess using $\MSO$ or, equivalently, tree-automata. 

We prove that we can check the correctness of such a guess by
introducing a special type of run, called \emph{loop}, which is
basically a run that starts and ends with the same stack. A run from one
milestone to the next will mainly consist of loops combined with a
finite number of stack operations.

\subsection{Milestones}

\begin{definition}[Milestone]
  A substack $s'$ of  $s=w_1:w_2:\ldots:w_n$ is a \emph{milestone} if
  $s'=w_1:w_2:\ldots :w_i:w'$ such that $0\leq i <n$ and $w_i\sqcap
  w_{i+1}\leq w'\leq w_{i+1}$.
  We denote by $\Milestones(s)$ the set of milestones of $s$. 
\end{definition}
Note that the substack relation $\leq$ linearly orders
$\Milestones(s)$.

\begin{lemma} \label{STACS:Lemma_VisitOrder} 
  If $s, t, m$ are stacks with $m\in \Milestones(t)$
    but $m\not\leq s$, then every run from $s$ to $t$ visits $m$. 
    Thus, for every run $r$ from the initial configuration to $s$, the
  function
  \begin{align*}
    f:\Milestones(s) \rightarrow \domain(r), & &
    s'\mapsto
    \max\{i\in\domain(r): r(i)=(q,s')\text{ for some }q\in Q\}    
  \end{align*}
  is an order
  embedding with respect to substack relation on the milestones and
  the natural order of $\domain(r)$.   
\end{lemma}   
In order to state the close correspondence between milestones of
a stack $s$ and the elements of $\Encode(s)$, we need the following definition.

\begin{definition}\label{STACS:DefLeftTree}
  Let $T\in\EncTrees$ be a tree and $d\in
  T\setminus\{\varepsilon\}$. Then the \emph{left and
  downward closed tree induced by $d$} is
  $\LeftTree{d,T} \coloneqq T{\restriction}_D$ where 
  \mbox{$D\coloneqq\{d'\in T: d'\leq_{lex} d\}\setminus\{\varepsilon\}$}. 
  Then we denote by
  $\LeftStack(d,T)\coloneqq\Decode(\LeftTree{d,T})$ the 
  \emph{ left stack induced by $d$}. 
\end{definition}
\begin{remark}
  $\LeftStack(d,s)$ is a substack of $s$ for all
  $d\in\domain(\Encode(s))$.
  This observation follows from Remark \ref{STACS:milestonesInEncoding}
  combined with the fact that the left stack is induced by a  
  lexicographically downward closed subset. In fact, $\LeftStack(d,s)$
  is a milestone of $s$.
\end{remark}

\begin{lemma} \label{STACS:LemmaOrderIso}
  The map given by $g:
  d \mapsto \LeftStack(d,\Encode(s))$ is an order isomorphism between 
  $\left(\domain(\Encode(q,s))\setminus\{\varepsilon\},\leq_{lex}\right)$
  and $\left(\Milestones(s),\leq\right)$.   
\end{lemma}  
Lemmas  \ref{STACS:LemmaOrderIso} and \ref{STACS:Lemma_VisitOrder} 
imply that every run  
$r$ decomposes as $r= r_1 \circ  r_2 \circ
\ldots \circ r_n $ where $r_i$ is a run from
the $i$-th milestone of $r(\length(r))$ to
the $(i+1)$-st milestone.  

In order to describe the structure of the $r_i$, we have to introduce
the notion of a loop. Informally speaking, a loop is a run $r$ that starts
and ends with the same stack $s$ and which does not look too much into $s$.

\begin{definition}
  Let $r$ be a run  of length $n$ with
  $r(i)=(q_i,s_i)$ for all $0\leq i \leq n$.
  \begin{itemize}
  \item 
    $r$ is called a \emph{simple high loop} if  
    $s_0=s_n$ and if $s_0<s_i$ for all $0<i<n$. 
  \item
    $r$ is called a \emph{simple low loop} of $s$ if
  $s_0=s_n=s$, between $0$ and $n$ the stack $s$ is
  never visited, $s_1 = \Pop{1}(s)$, $\Lvl(s)=1$, 
  $\lvert s_i\rvert \geq \lvert s \rvert$ for all $0\leq i \leq n$,
  and $r{\restriction_{[2,n-1]}}$ is the composition of
  simple low loops and simple high loops of $\Pop{1}(s)$.
\item
  $r$ is called \emph{loop} if it is a finite composition of low loops
  and high loops.
  \end{itemize}
\end{definition}

\begin{lemma}\label{STACS:ClonePopMilestoneRun}
  Let $s$ be some stack, $m_1, m_2$  milestones of $s$, and $r$  a
  run from $m_1$ to $m_2$ that never visits
  any other milestone of $s$. 
  Then either $r= l_1 \circ p\circ l_2$  or
  $r= l_0 \circ c \circ l_1 \circ p_1 \circ l_2 \circ p_2 \circ l_3
  \circ \ldots \circ p_{n} \circ l_{n+1}$ where 
  each $l_i$ is a loop, and all $p_i,
  p$, and $c$ are runs of length $1$, $p$ performs one
  $\Push{\sigma,k}$, $c$ performs one $\Clone{2}$, and the $p_i$
  perform one $\Pop{1}$ each. 
\end{lemma}

This lemma motivates why we only define low loops for stacks $s$ with
$\Lvl(s)=1$. Whenever the topmost symbol of a milestone $m$ is not a
cloned element, then $\Pop{1}(m)$ is another milestone. 
Hence, the $l_i$ can only contain low loops if they start at a stack
with cloned topmost symbol. But any stack $s$ with cloned topmost
symbol and $\Lvl(s)=2$ cannot be restored from $\Pop{1}(s)$ without
passing $\Pop{2}(s)$ since a $\Push{\sigma,2}$-operation
would create the wrong link-level.    

From Lemma \ref{STACS:ClonePopMilestoneRun} we can derive 
that  deciding whether there is a run
from one milestone to the next is possible
if we know the pairs of initial and final states of loops
of certain stacks $s$. Hence we are interested in the sets
$\Loops(s) \subseteq Q\times Q$ with
$(q_1,q_2)\in \Loops(s)$ if  and only if there is a loop from 
$(q_1,s)$  to $(q_2,s)$.
The crucial observation is that $\Loops(s)$ may be
calculated by a finite automaton reading $\TOP{2}(s)$. 

\begin{lemma}
  \label{STACS:DetermineLoops}
  For every \CPS there exists a finite automaton $A$ that
  calculates\footnote{ We consider the final state reached by $A$
    on input $w$ as the value it calculates for $w$.}
  on input \mbox{$w\in (\Sigma \times\{1,2\})^*$} the set $\Loops(s)$ for all
  stacks $s$ such that $w=\pi(\TOP{2}(s))$. Here,  
  \mbox{$\pi:(\Sigma\times\{1,2\}\times\N)^*\rightarrow
    (\Sigma\times\{1,2\})^*$} is 
  the projection onto the symbols and collapse-levels.
\end{lemma}

\subsection{Detection of Reachable Configurations}

We have already seen that every run to a valid configuration $(q,s)$
passes all the milestones of $s$. Now, we use the last state in 
which a run $r$ to $(q,s)$ visits each milestone as a certificate for
the reachability of $(q,s)$. To be precise, \emph{a certificate for the
reachability of $(q,s)$}  is a map
$f:\domain\big(\Encode(q,s)\big)\setminus\{\varepsilon\} \rightarrow Q$ such
that there is some run $r$ from $\bot_2$ to $(q,s)$ and $f(d)=q$ if
and only if $r(i)=\big(q,\LeftStack(d)\big)$ for $i$ 
the maximal position in $r$ where $\LeftStack(d)$ is visited.

\begin{lemma}\label{STACS:CertificateCheck}
  For every \CPG $G$, there is a tree-automaton
  that checks for each map
  \begin{align*}
   f:\domain(\Encode(q,s))\setminus\{\varepsilon\} &\rightarrow Q 
  \end{align*}
  whether $f$ is a certificate of the reachability of $(q,s)$, i.e.,
  whether $f$ is induced by some
  run $r$ from the initial configuration to $(q,s)$. 
\end{lemma}

The proof of the lemma uses Lemma \ref{STACS:DetermineLoops} and the fact
that the path from the root to some $d\in\Encode(s)$ encodes the
topmost word of $\LeftStack(d,\Encode(s))$. Hence, a tree automaton
reading $\Encode(s)$ is able to calculate for each position
$d\in\Encode(s)$ the pairs of initial and final states of loops of
$\LeftStack(d)$. As every run decomposes as a sequence of loops
separated by a single operation, knowing $\Loops(s')$ for each $s'\leq
s$ enables the automaton to check the correctness of a candidate for a
certificate of reachability. 

As a tree-automaton may non-deterministically guess a
certificate of the reachability of a configuration, 
the encodings of reachable configurations form a
regular set.

\subsection{Extension to Regular Reachability} 
By now, we have already established the tree-automaticity of each
\CPG $G$ since we have seen that our encoding yields a regular image of
the vertices of $G$ and the transition relations are turned into
regular relations of the tree encoding. 
Using similar techniques, we can improve this result: 

\begin{theorem} \label{STACS:CPGwithREACHTreeAutomatic}
  If $G$ is the $\varepsilon$-closure
  of some
  $\CPG$ $G'$ then $(G, \Reach)$ is tree-automatic where $\Reach$ is
  the binary predicate that is true on a pair $(c_1,c_2)$ of
  configurations if there is a path from $c_1$ to $c_2$ in $G$. 
\end{theorem}
\begin{remark}
  Each graph in the second level of the Caucal-hierarchy can be
  obtained as the $\varepsilon$-contraction of some level $2$
  \CPG (see \cite{cawo03}) whence all these graphs are tree-automatic.
\end{remark}

For a \CPS $S$ let $R\subseteq\Delta^*$ be a regular language
over the transitions of $S$. As collapsible pushdown graphs are closed
under products with finite automata even the reachability predicate
$\Reach_R$ with restriction to $R$ is
tree-automatic. Here, $\Reach_R xy$ holds if there is a path from $x$
to $y$ in $\Graph(S)$ that uses a sequence of transitions in $R$. If $A$
is the automaton recognising $R$, we obtain that $\Reach_R
(q,s)(q',s')$ holds in $\Graph(S)$ iff $\Reach
\big((q,q_i),s\big)\big((q',q_f),s'\big)$ holds
in $\Graph(S\times A)$ where 
$q_i$ is the initial and $q_f$ the unique final state of
$A$. Using this idea one can define a \CPG $G'$ which is basically
$\Graph(S\cup (S\times A))$ extended by transitions from $(q,s)$ to
$((q,q_i),s)$ and to $((q,q_f),s)$. $\Graph(S)$ as well as
$\Reach_R$ w.r.t. $\Graph(S)$ are $\FO[\Reach]$-interpretable in $G'$. Hence
we obtain:

\begin{theorem}
  Given a collapsible pushdown graph of level $2$, its $\FO[\Reach_R]$
  theory is decidable for each regular $R\subseteq \Delta^*$.  
\end{theorem}

\subsection{Computation of concrete tree-automatic representations of
CPG}

Up to now, we have only seen that there is a tree-automatic
representation for each \CPG. For computing a concrete representation,
we rely on the following lemma.

\begin{lemma} \label{STACS:ReachableConfigsinMuCalcul}
  Given some  \CPS
  \mbox{$S=(\Gamma, Q, \Delta, q_0)$}, some $q\in Q$, and some 
  stack $s$,
  it is decidable whether $(q,s)$ is a vertex of $\Graph(S)$.
\end{lemma}
The proof is based on the idea that a stack is uniquely determined by
its top element and the information which substacks can be reached via
$\Collapse$- and $\Pop{i}$-operations. Hence we can construct an
extension $S'$ of $S$ and a modal formula $\varphi_{q,s}$ such that there
is some element $v\in \Graph(S')$ satisfying
$\Graph(S'),v \models\varphi_{q,s}$ iff $(q,s)\in\Graph(S)$. $S'$
basically contains new states for every substack of $s$ and connects
the different states via the appropriate $\Pop{i}$-operations which
are only applied if the topmost symbol of the stack agrees with the
symbol we would expect when starting the $\Pop{i}$-sequence in
configuration $(q,s)$. 

From this lemma we can derive the computability of the automata in
Lemma \ref{STACS:DetermineLoops}. Having obtained these automata, 
the construction of a tree-automatic representation of some \CPG is
directly derived from the proofs yielding the following theorem.

\begin{theorem}
  There is an algorithm that, given a level $2$
  \CPG $G$ and regular sets \mbox{$R_1, \ldots, R_n
    \subseteq\Delta^*$}, 
   computes a tree-automatic representation of
  $(G,\Reach_{R_1},\ldots,\Reach_{R_n})$. 
\end{theorem}

\section{Conclusion}
We have seen that level $2$ collapsible pushdown graphs are
tree-automatic. This result holds also if we apply
$\varepsilon$-contractions and if we add regular reachability
predicates.
This implies that the second level of the Caucal-hierarchy is
tree-automatic.  
But our result can only be seen as a starting point for further
investigations of the \CPG hierarchy:
are level $3$
collapsible pushdown graphs tree-automatic? 
We know an example of a level $5$ \CPG which is not tree-automatic. 
But even when tree-automaticity of all \CPG cannot be expected, the
question remains whether all \CPG have decidable
\FO theories. 
In order to solve this problem one has to come up with 
new techniques. 

A rather general question concerning our result aims at our knowledge
about tree-automatic structures. Recent developments 
in the string case \cite{Kuske2009} show the decidability of rather large
extensions of first-order logic for automatic structures. 
It would be interesting to clarify the status of the analogous claims
for tree-automatic structures. 
Positive answers concerning the
decidability of extensions of first-order logic on tree-automatic
structures would give us the
corresponding decidability results for collapsible pushdown graphs of
level $2$. 

\bibliography{KartzowBibStacs10}
\bibliographystyle{plain}

\end{document}